\documentclass[12pt,a4paper]{amsart}
\usepackage{amsfonts}
\usepackage{color}
\usepackage{amsthm}
\usepackage{amsmath}
\usepackage{amscd}
\usepackage[latin2]{inputenc}
\usepackage{t1enc}
\usepackage[mathscr]{eucal}
\usepackage{indentfirst}
\usepackage{graphicx}
\usepackage{graphics}
\usepackage{pict2e}
\usepackage{epic}
\usepackage[margin=2.9cm]{geometry}
\usepackage{epstopdf}
\usepackage{amssymb,amsthm, graphicx, stmaryrd, fancyhdr}
\usepackage{amsfonts}
\usepackage{amsmath}
\usepackage{epsfig,epstopdf}
\usepackage{subfigure}
\usepackage{latexsym}
\usepackage{color}

\numberwithin{equation}{section}

\theoremstyle{plain}

\newtheorem{theorem}{Theorem}[section]
\newtheorem{lemma}{Lemma} [section]

\theoremstyle{definition}

\newtheorem{remark}{Remark}[section]

\begin{document}
\title{A nonlinear population model}
\author{Dragos-Patru Covei$^1$}
\author{Traian A. Pirvu$^2$}
\author{Catalin Sterbeti$^3$}
\date{\\
$^1$Department of Applied Mathematics, The Bucharest University of Economic
Studies, Piata Romana, 1st district, postal code: 010374, postal office: 22,
Romania\\
$^2$Department of Mathematics and Statistics, McMaster University, 1280 Main
Street West, Hamilton, ON, L8S 4K1, Canada \\
$^3$ Department of Applied Mathematics, University of Craiova, 13, A.I. Cuza
Street, 200585, Craiova, Dolj, Romania. }
\email{coveidragos@yahoo.com; tpirvu@math.mcmaster.ca; sterbetiro@yahoo.com}

\begin{abstract}
This paper considers a nonlinear model for population dynamics with age
structure. The fertility rate with respect to age is non constant and has
the form proposed by \cite{1}. Moreover, its multiplicative structure and
the multiplicative structure of mortality makes the model separable. In this
setting it is shown that the number of births in unit time is given by a
system of nonlinear ordinary differential equations. The steady state
solution together with the equilibrium solution is found explicitly.
\end{abstract}

\maketitle

\section{Introduction \label{intro}}

Our paper studies a population model with age structure. This model and its
variants were considered in many works such as \cite{Bra1, Bra2, Hop, GG1,
Ian, McKen, ShaLot, Ster,W}, to cite just a few. There is a vast literature
on this subject by now.

The seminal work of \cite{McKen} considers the age structure into the
dynamics of one sex population model assuming that the female population
dynamics can be modelled as a function of two variables, namely age and
time. This model takes as inputs an age specific mortality intensity and an
age specific fertility function. The number of individuals at a given time
that have age less than a certain age is given by an integral of a function $%
p(a,t)$ of two variables (age and time). The size of population (the total
number of individuals) can also be obtained by integrating $p(a,t).$ It
appears only natural to find $p(a,t)$ which is described by a system of
integral and differential equations. This system is shown to reduce to a
Volterra integral equation. The later has a unique solution established by a
fixed point theorem, but this solution is not available in closed form, and
its numerical computation based on fixed point iterations seems complicated.
The next class of population models are the nonlinear models obtained in
situations when the mortality and fertility are functions of age and
population size (see \cite{2} for more on this).

Let us mention two seminal works in the paradigm of nonlinear population
models. The first work to incorporate population size dependent mortality
and fertility rates (rendering the model nonlinear) is \cite{GG1}. They also
characterize the equilibrium population and established a condition for the
equilibrium to be locally asymptotically stable. A special class of models,
the separable ones, are considered in their paper. \cite{W} shows that the
stability classification, depends in many cases on the marginal birth and
death rates as measures the sensitivities of the fertility and mortality; on
other cases more information is required to determine the stability.

Let us turn now to the contributions of our paper. Our setting is the
nonlinear population model, in the special case when mortality and fertility
are separable functions of age and population size. The age part of the
fertility is assumed as in \cite{1} which is a model with non constant (in
age) fertility rate. In this paradigm we establish the existence of the
steady state for the nonlinear equations characterizing the population
dynamics. Moreover, these are shown to be equivalent to a nonlinear system
of differential equations. The late is shown to have an equilibrium solution
which we find explicitly. The equilibrium stability analysis can also be
established.

The reminder of this paper is organized as follows: Section 2 presents the
populations models, linear, nonlinear. Sections 3 presents our nonlinear
model and the main results of the paper.

\section{Population Models}

Let us first introduce the linear population model. The exposition here
follows \cite{2}. The dynamics of population is expressed in terms of the
density of the population of age $a$, at time $t$ density denoted $p(a,t)$.
The total population at $t$, denoted by $P(t)$ can be obtained by
integrating its density, i.e., 
\begin{equation*}
P(t)=\int_{0}^{\infty }p(a,t)da.
\end{equation*}%
Next let us introduce \textbf{fertility} and \textbf{mortality}. The number
of offspring, borne by individuals during the infinitesimal time-interval $%
[t,t+dt]$, and the infinitesimal age-interval $[a,a+da]$ is 
\begin{equation*}
\beta (a,t),
\end{equation*}%
also referred to as the age specific fertility. The number of offspring
during the infinitesimal interval $[t,t+dt]$ is then 
\begin{equation*}
B(t)=\int_{0}^{\infty }\beta (a,t)p(a,t)da.
\end{equation*}%
The number of deaths of individuals during the infinitesimal time-interval $%
[t,t+dt]$, and the infinitesimal age-interval $[a,a+da]$ is 
\begin{equation*}
\mu (a,t).
\end{equation*}%
The number of deaths during the infinitesimal interval $[t,t+dt]$ is then 
\begin{equation*}
D(t)=\int_{0}^{\infty }\mu (a,t)p(a,t)da.
\end{equation*}%
The probability that an individual of age $a-x$ at the time $t-x$ will
survive up to time $t$ (with age $a$) is given by 
\begin{equation*}
\pi (a,t,x)=e^{-\int_{0}^{x}\mu (a-\sigma ,t-\sigma )d\sigma }.
\end{equation*}%
In the case of time independent mortality 
\begin{equation*}
\pi (a)=e^{-\int_{0}^{a}\mu (\sigma )d\sigma },
\end{equation*}%
is the probability for a newborn to survive to age $a$, also known as the
survival probability.

In the following we will derive the linear Lotka-McKendrick Equation. The 
\textbf{fertility} and \textbf{mortality} rates $\beta(a),$ and $\mu(a)$ are
assume time independent, as they only depend on the age $a.$ The number of
individuals with age less than $a$ at time $t,$ denoted by $N(a,t)$ is given
by

\begin{equation*}
N(a,t)= \int_{0}^{a} p(\sigma,t) d\sigma.
\end{equation*}

Next, let us look at the number of individuals with age less than $a+h$ at
time $t+h,$ i.e., $N(a+h,t+h)$. This number will comprise $N(a,t),$ and all
newborn in the time interval $[t,t+h]$ (their age will be less than $a+h$),
which is 
\begin{equation*}
\int_{t}^{t+h}B(s)ds.
\end{equation*}%
One needs to adjust then for the deaths from the newborns, through the time
interval $[t,t+h],$ and the deaths on $[t,t+h]$ of individuals older than $%
a, $ and these number of deaths is 
\begin{equation*}
\int_{0}^{h}\int_{0}^{a+s}\mu (\sigma )p(\sigma ,t+s)d\sigma ds.
\end{equation*}%
This is the case because 
\begin{equation*}
\int_{0}^{a+s}\mu (\sigma )p(\sigma ,t+s)d\sigma ,
\end{equation*}%
gives the number of individuals who die at $t+s$ younger than $a+s$.
Therefore 
\begin{equation*}
N(a+h,t+h)=N(a,t)+\int_{t}^{t+h}B(s)ds-\int_{0}^{h}\int_{0}^{a+s}\mu (\sigma
)p(\sigma ,t+s)d\sigma ds.
\end{equation*}%
Let us differentiate this with respect to $h$, and then take $h=0$ to get 
\begin{equation*}
p(a,t)+\int_{0}^{a}p_{t}(\sigma ,t)d\sigma =B(t)-\int_{0}^{a}\mu (\sigma
)p(\sigma ,t)d\sigma .
\end{equation*}%
Next, let us differentiate this with respect to $a,$ to get 
\begin{equation*}
p(a,t)+p_{a}(a,t)+\mu (a)p(a,t)=0.
\end{equation*}%
Also by setting $a=0$ yields 
\begin{equation*}
p(0,t)=B(t),
\end{equation*}%
but on the other hand 
\begin{equation*}
B(t)=\int_{0}^{\infty }\beta (\sigma )p(\sigma ,t)d\sigma .
\end{equation*}%
As such we obtained the following system 
\begin{equation}
\left\{ 
\begin{array}{lll}
p_{a}(a,t)+p_{t}(a,t)+\mu (a)p(a,t)=0 & \text{for } & a,t\geq 0, \\ 
p(0,t)=\int\limits_{0}^{\infty }\beta (\sigma )\rho (\sigma ,t)\,d\sigma & 
\text{for } & t>0, \\ 
p(a,0)=p_{0}(a) & \text{for} & a>0.%
\end{array}%
\right.  \label{ddts}
\end{equation}

\subsection{A Special Linear Population Model}

Inspired by \cite{NPED}, \cite{Dufresne} and \cite{F}, the recent work \cite%
{CPS} considers a model with survival probability $\pi (a)$ pseudo
exponential, i.e. 
\begin{equation}
\pi (a)=\sum_{i=1}^{n}c_{i}e^{{-\mu _{i}a}},  \label{s1}
\end{equation}%
for positive constants $c_{i},\mu _{i}$ with 
\begin{equation*}
\sum_{i=1}^{n}c_{i}^{2}\neq 0\text{ ,}\sum_{i=1}^{n}\mu _{i}^{2}\neq 0\text{
and }\sum_{i=1}^{n}c_{i}=1.
\end{equation*}%
Moreover, we present the case of 
\begin{equation}
\pi (a)=\left( \sum_{i=0}^{n}c_{i}a^{i}\right) e^{{-\mu _{1}a}}.  \label{s2}
\end{equation}%
In this setting finding $p(a,t)$ is reduced to a linear ODE system. In
special situations\ (\ref{s1})-(\ref{s2}) a closed form solution is obtained
by means of Laplace transform.

\subsection{Nonlinear Population Models}

These are models in which death rate $\mu (a,p),$ and the fertility rate $%
\beta (a,p),$ are functions of age and population size. The age density
function at time $t\geq 0,$ $p\left( a,t\right) $ satisfies the following
nonlinear equations

\begin{equation*}
\left\{ 
\begin{array}{l}
p_{t}\left( a,t\right) +p_{a}\left( a,t\right) +\mu (a,P(t))p\left(
a,t\right) =0 \\ 
p\left( 0,t\right) =\int_{0}^{\infty }\beta (\sigma ,P(t))p\left( \sigma
,t\right) d\sigma \\ 
p\left( a,0\right) =p_{0}\left( a\right) \\ 
P\left( t\right) =\int_{0}^{\infty }p\left( \sigma ,t\right) d\sigma ,%
\end{array}%
\right.
\end{equation*}%
where $P(t)$ is the total population at $t$. This system of nonlinear
equations can be reduced to integral equations, as it is shown in (see \cite%
{2} for more details). Indeed, let 
\begin{equation*}
\pi (a,t,x,P)=e^{-\int_{0}^{x}\mu (a-\sigma ,P(t-\sigma ))d\sigma }
\end{equation*}%
and

\begin{equation}
p(a,t)=\left\{ 
\begin{array}{lcr}
p_{0}(t-a)\pi (a,t,t,P) & \mbox{ for } & t\leq a \\ 
B(t-a)\pi (a,t,a,P) & \mbox{ for
} & t>a.%
\end{array}%
\right.  \label{iss}
\end{equation}%
Here the function $B$ solves the following system of equations:

\begin{equation}
\left\{ 
\begin{array}{lcr}
B(t)=\int_{0}^{t}\beta (\sigma ,P(t))\pi (\sigma ,t,\sigma ,P)B(\sigma
)d\sigma +F(t,P) &  &  \\ 
P(t)=\int_{0}^{t}\pi (\sigma ,t,\sigma ,P)B(\sigma )d\sigma +G(t,P) &  & 
\end{array}%
\right.  \label{iss1}
\end{equation}%
where

\begin{equation}
\left\{ 
\begin{array}{lcr}
F(t,P)=\int_{t}^{\infty }\beta (\sigma ,P(t))\pi (\sigma ,t,t,P)p_{0}(\sigma
-t)d\sigma &  &  \\ 
G(t,P)=\int_{t}^{\infty }\pi (\sigma ,t,t,P)p_{0}(\sigma -t)d\sigma . &  & 
\end{array}%
\right.  \label{iss2}
\end{equation}%
This system (\ref{iss1}) can be solved through the following iterative method

\begin{equation}
\left\{ 
\begin{array}{lcr}
B^{k+1}(t)=\int_{0}^{t}\beta (\sigma ,P^{k}(t))\pi (\sigma ,t,\sigma
,P^{k})B^{k}(\sigma )d\sigma +F(t,P^{k}) &  &  \\ 
P^{k+1}(t)=\int_{0}^{t}\pi (\sigma ,t,\sigma ,P^{k})B^{k}(\sigma )d\sigma
+G(t,P^{k}), &  & 
\end{array}%
\right.  \label{iss11}
\end{equation}%
see \cite{2} for more details on this. A special class of nonlinear models
are the separable population models presented in the next subsection.

\subsection{Separable Population Models}

Now let us specialize the nonlinear model with the following choice of
fertility and mortality 
\begin{equation*}
\beta (a,p)=R_{0}\beta _{0}\Phi \left( p\right) e^{-\rho a},\quad \mu
(a,p)=\mu _{0}+\Psi (p).
\end{equation*}%
By plugging this in the nonlinear system one gets

\begin{equation*}
\left\{ 
\begin{array}{l}
p_{t}\left( a,t\right) +p_{a}\left( a,t\right) +(\mu _{0}+\Psi (P(t)))\left(
a,t\right) =0 \\ 
p\left( 0,t\right) =\int_{0}^{\infty }R_{0}\Phi \left( P(t)\right) e^{-\rho
a}p\left( \sigma ,t\right) d\sigma \\ 
p\left( a,0\right) =p_{0}\left( a\right) \\ 
P\left( t\right) =\int_{0}^{\infty }p\left( \sigma ,t\right) d\sigma .%
\end{array}%
\right.
\end{equation*}%
These nonlinear equations can be reduced to the following ODE system

\begin{equation}
\left\{ 
\begin{array}{l}
P^{\prime }\left( t\right) =-\left( \mu _{0}+\Psi \left( P\left( t\right)
\right) \right) P\left( t\right) +R_{0}\beta _{0}\Phi \left( P\left(
t\right) \right) Q\left( t\right) \\ 
Q^{\prime }\left( t\right) =\left( R_{0}\beta _{0}\Phi \left( P\left(
t\right) \right) -\rho -\mu _{0}-\Psi \left( P\left( t\right) \right)
\right) Q\left( t\right)%
\end{array}%
\right.  \label{stud1}
\end{equation}%
where 
\begin{equation*}
Q\left( t\right) =\int_{0}^{\infty }e^{-a\rho }p_{t}\left( a,t\right) da,
\end{equation*}%
see \cite{2} for more details.

Let us turn now to the existence of steady states. The net reproduction
number is given by 
\begin{equation*}
R(x)=R_{0}\Phi (x)\int_{0}^{\infty }\beta _{0}e^{-\left( \rho +\mu _{0}+\Psi
\left( x\right) \right) a}da.
\end{equation*}%
This quantity was introduced first by \cite{GG1}. According with \cite{GG1}
the quantity $R\left( x\right) $ is the number of children expected to be
born to an individual when the population is $x$. The steady state $P^{\ast
} $ is given by the following equation 
\begin{equation*}
R(P^{\ast })=1,
\end{equation*}%
which in this setting reads

\begin{equation*}
R_0 \beta_0 \Phi(P^{*})= \rho +\mu _{0}+\Psi\left( P^{*}\right).
\end{equation*}

\section{A Special Nonlinear Population Model}

We consider the following logistic system%
\begin{equation}
\left\{ 
\begin{array}{l}
p_{t}\left( a,t\right) +p_{a}\left( a,t\right) +\mu _{0}p\left( a,t\right)
+\Psi \left( P\left( t\right) \right) p\left( a,t\right) =0 \\ 
p\left( 0,t\right) =R_{0}\Phi \left( P\left( t\right) \right)
\int_{0}^{\infty }\overset{n-1}{\underset{i=0}{\Sigma }}\beta _{i}\sigma
^{i}e^{-\rho \sigma }p\left( \sigma ,t\right) d\sigma \\ 
p\left( a,0\right) =p_{0}\left( a\right) \\ 
P\left( t\right) =\int_{0}^{\infty }p\left( \sigma ,t\right) d\sigma%
\end{array}%
\right.  \label{SNPM}
\end{equation}%
where\ $p\left( a,t\right) $ is the age density function at time $t\geq 0$, $%
P\left( t\right) $ is the total population at time $t$, $p_{0}\left(
a\right) $ is the initial population at time $t=0$, $a\in \left[ 0,\infty
\right) $, $\mu _{0}>0$ is the intrinsic mortality term and $R_{0}$, $\rho $%
, $\sigma $ are prescribed positive parameters. In real world, the age
profile is given by%
\begin{equation*}
\frac{p\left( a,t\right) }{P\left( t\right) }.
\end{equation*}%
The fertility is given by 
\begin{equation*}
\beta \left( a,p\right) =R_{0}\beta _{0}\Phi \left( p\right) \overset{n-1}{%
\underset{i=0}{\Sigma }}\beta _{i}a^{i}e^{-\rho a}\text{ with }\beta
_{0},\beta _{1},...,\beta _{n-1}\in \left( 0,\infty \right) ,
\end{equation*}%
and the mortality by 
\begin{equation*}
\mu (a,P)=\mu _{0}+\Psi (P).
\end{equation*}%
The couple\ $\left( \beta \left( a\right) ,\mu _{0}\right) $ can be
interpreted as an intrinsic birth-death process that is age dependent while $%
\Psi \left( P\left( t\right) \right) $ models an external mortality that is
the same for all ages and just depends on the weighted sizes.

Since we are interested in the existence of steady states solutions for the
system (\ref{SNPM}), i.e. for solutions that are constant in time, we assume
for the start that\ $\Phi $ and $\Psi $ are continuous on $\left[ 0,\infty
\right) $ continuously differentiable on $\left( 0,\infty \right) $ and%
\begin{eqnarray}
\Phi \left( x\right) &\geq &0,\text{ }\Phi ^{\prime }\left( x\right) <0\text{%
, }\Phi \left( 0\right) =1,\Phi \left( +\infty \right) =0,  \label{A1} \\
\Psi \left( x\right) &\geq &0,\Psi ^{\prime }\left( x\right) >0,\Psi \left(
0\right) =0,\Psi \left( +\infty \right) =+\infty .  \label{A2}
\end{eqnarray}%
We also adopt the normalization condition 
\begin{equation*}
\int_{0}^{\infty }\overset{n-1}{\underset{i=0}{\Sigma }}\beta
_{i}a^{i}e^{-\left( \rho +\mu _{0}\right) a}da=1,
\end{equation*}%
from where, for example with the use of Gamma integral, we obtain%
\begin{equation}
\overset{n-1}{\underset{i=0}{\Sigma }}\frac{\beta _{i}i!}{\left( \rho +\mu
_{0}\right) ^{i+1}}=1,  \label{normc}
\end{equation}%
so that the parameter $R_{0}$ in (\ref{SNPM}) has the role of an intrinsic
basic reproduction number denoted by $R\left( x\right) $ in the next.

Concerning the existence of steady states, the net reproduction number at
size $x$ (see page 154 in \cite{2} or page 288 in \cite{GG1}) takes the
following form within our setting 
\begin{eqnarray*}
R\left( x\right) &=&R_{0}\Phi \left( x\right) \int_{0}^{\infty }\overset{n-1}%
{\underset{i=0}{\Sigma }}\beta _{i}a^{i}e^{-\left( \rho +\mu _{0}+\Psi
\left( x\right) \right) a}da \\
&=&R_{0}\Phi \left( x\right) \overset{n-1}{\underset{i=0}{\Sigma }}\beta
_{i}\int_{0}^{\infty }a^{i}e^{-\left( \rho +\mu _{0}+\Psi \left( x\right)
\right) a}da \\
&=&R_{0}\Phi \left( x\right) \overset{n-1}{\underset{i=0}{\Sigma }}\beta
_{i}\int_{0}^{\infty }\left( \frac{t}{\rho +\mu _{0}+\Psi \left( x\right) }%
\right) ^{i}e^{-t}\frac{dt}{\rho +\mu _{0}+\Psi \left( x\right) } \\
&=&R_{0}\Phi \left( x\right) \overset{n-1}{\underset{i=0}{\Sigma }}\frac{%
\beta _{i}}{\left( \rho +\mu _{0}+\Psi \left( x\right) \right) ^{i+1}}%
\int_{0}^{\infty }t^{i}e^{-t}dt \\
&=&R_{0}\Phi \left( x\right) \overset{n-1}{\underset{i=0}{\Sigma }}\frac{%
\beta _{i}\Gamma \left( i+1\right) }{\left( \rho +\mu _{0}+\Psi \left(
x\right) \right) ^{i+1}} \\
&=&R_{0}\Phi \left( x\right) \overset{n-1}{\underset{i=0}{\Sigma }}\frac{%
\beta _{i}i!}{\left( \rho +\mu _{0}+\Psi \left( x\right) \right) ^{i+1}}
\end{eqnarray*}%
where $\Gamma $ represent the Gamma integral. A first observation regarding
the set of assumptions (\ref{A1})-(\ref{A2}) is that 
\begin{equation}
\lim_{x\rightarrow \infty }R\left( x\right) =0\text{ and }R\left( x\right) 
\text{ is a decreasing function (i.e. }R^{\prime }\left( x\right) <0\text{).}
\label{zdec}
\end{equation}%
Indeed $R\left( x\right) $ is a decreasing function because 
\begin{equation*}
R^{\prime }(x)=\overset{n-1}{\underset{i=0}{\Sigma }}\frac{\beta _{i}i!\left[
R_{0}\Phi ^{\prime }\left( x\right) \left( \rho +\mu _{0}+\Psi \left(
x\right) \right) -R_{0}\Phi \left( x\right) \Psi ^{\prime }\left( x\right) %
\right] }{\left( \rho +\mu _{0}+\Psi \left( x\right) \right) ^{i+2}}-\overset%
{n-1}{\underset{i=0}{\Sigma }}\frac{i\beta _{i}i!R_{0}\Phi \left( x\right)
\Psi ^{\prime }\left( x\right) }{\left( \rho +\mu _{0}+\Psi \left( x\right)
\right) ^{i+2}}<0\text{,}
\end{equation*}%
for all $x\geq 0$. Moreover, 
\begin{equation*}
\lim_{x\rightarrow \infty }R\left( x\right) =0,
\end{equation*}%
is satisfied in light of the asymptotic conditions on $\Psi ,\Phi .$

As is well known (see page 154 in \cite{2} or \cite[Theorem 6, pages 288-289]%
{GG1}), since we have a single weighted size $P$, we use that non-trivial
stationary sizes $P^{\ast }$ must satisfy

\begin{equation}
R\left( P^{\ast }\right) =1,  \label{nsc}
\end{equation}%
and this is a necessary and sufficient condition for a non-trivial
stationary sizes to exist with total population $P^{\ast }$. In this case, (%
\ref{nsc}) becomes%
\begin{equation}
\frac{R_{0}\Phi \left( P^{\ast }\right) }{\rho +\mu _{0}+\Psi \left( P^{\ast
}\right) }\overset{n-1}{\underset{i=0}{\Sigma }}\frac{\beta _{i}i!}{\left(
\rho +\mu _{0}+\Psi \left( P^{\ast }\right) \right) ^{i}}=1.  \label{ec}
\end{equation}
Due to (\ref{zdec}) this equation (\ref{ec}) has one, and only one,
non-trivial solution if, and only if, $R_{0}>1$. The fact that the condition 
$R_{0}>1$ is necessary and sufficient for the existence of a non-trivial
equilibrium means that $R_{0}$ acts as a bifurcation parameter: under the
assumptions of the model it is clear that we have a forward bifurcation at
the point $R_{0}$.

One can rewrite (\ref{ec}) for the non-trivial stationary sizes $P^{\ast }$
in the form%
\begin{equation}
\overset{n-1}{\underset{i=0}{\Sigma }}\frac{\beta _{i}i!}{\left( \rho +\mu
_{0}+\Psi \left( P^{\ast }\right) \right) ^{i}}=\frac{\rho +\mu _{0}+\Psi
\left( P^{\ast }\right) }{R_{0}\Phi \left( P^{\ast }\right) } .
\label{verif}
\end{equation}

Let us summarize the results in the following Lemma.

\begin{lemma}
Given assumptions \eqref{A1}, \eqref{A2} our model has a unique steady state 
$P^*$ given by equation \eqref{ec} if and only if $R_0>1.$
\end{lemma}

Now let us turn to the problem of finding $P(t).$ We will denote by%
\begin{equation*}
P_{i}\left( t\right) =\int_{0}^{\infty }\sigma ^{i-1}e^{-\rho \sigma
}p\left( \sigma ,t\right) d\sigma ,
\end{equation*}%
for $i=1,2,..n$. The next step, is to observe that the \textit{renewal
condition }or\textit{\ the total birth rate }or\textit{\ fertility rate, }at
the time\textit{\ }$t$ can be can be written in the new notations such 
\begin{eqnarray*}
p\left( 0,t\right) &=&R_{0}\Phi \left( P\left( t\right) \right) \overset{n-1}%
{\underset{i=0}{\Sigma }}\beta _{i}\int_{0}^{\infty }\sigma ^{i}e^{-\rho
\sigma }p\left( \sigma ,t\right) d\sigma \\
&=&R_{0}\Phi \left( P\left( t\right) \right) \overset{n-1}{\underset{i=0}{%
\Sigma }}\beta _{i}P_{i+1}\left( t\right) .
\end{eqnarray*}%
More that, the calculation of the first derivative of $P\left( t\right) $
gives%
\begin{eqnarray*}
P^{\prime }\left( t\right) &=&\int_{0}^{\infty }p_{t}\left( a,t\right) da \\
&=&-\int_{0}^{\infty }p_{a}\left( a,t\right) da-\left( \mu _{0}+\Psi \left(
P\left( t\right) \right) \right) \int_{0}^{\infty }p\left( a,t\right) da \\
&=&p\left( 0,t\right) -\left( \mu _{0}+\Psi \left( P\left( t\right) \right)
\right) P\left( t\right) \\
&=&-\left( \mu _{0}+\Psi \left( P\left( t\right) \right) \right) P\left(
t\right) +R_{0}\Phi \left( P\left( t\right) \right) \overset{n-1}{\underset{%
i=0}{\Sigma }}\beta _{i}P_{i+1}\left( t\right)
\end{eqnarray*}%
and, similarly%
\begin{eqnarray*}
P_{1}^{\prime }\left( t\right) &=&\int_{0}^{\infty }e^{-a\rho }p_{t}\left(
a,t\right) da \\
&=&-\int_{0}^{\infty }e^{-a\rho }p_{a}\left( a,t\right) da-\left( \mu
_{0}+\Psi \left( P\left( t\right) \right) \right) \int_{0}^{\infty
}e^{-a\rho }p\left( a,t\right) da \\
&=&p\left( 0,t\right) -\rho \int_{0}^{\infty }e^{-a\rho }p_{a}\left(
a,t\right) da-\left( \mu _{0}+\Psi \left( P\left( t\right) \right) \right)
P_{1}\left( t\right) \\
&=&R_{0}\Phi \left( P\left( t\right) \right) \overset{n-1}{\underset{i=0}{%
\Sigma }}\beta _{i}P_{i+1}-\rho P_{1}\left( t\right) -\left( \mu _{0}+\Psi
\left( P\left( t\right) \right) \right) P_{1}\left( t\right) \\
&=&\left( R_{0}\beta _{0}\Phi \left( P\left( t\right) \right) -\rho -\mu
_{0}-\Psi \left( P\left( t\right) \right) \right) P_{1}\left( t\right)
+R_{0}\Phi \left( P\left( t\right) \right) \overset{n-1}{\underset{i=1}{%
\Sigma }}\beta _{i}P_{i+1}.
\end{eqnarray*}%
In the same way, the calculation of the first derivative of $P_{i+1}\left(
t\right) $ ($i=1,...,n-1$) gives 
\begin{eqnarray*}
P_{i+1}^{\prime }\left( t\right) &=&\int_{0}^{\infty }a^{i}e^{-\rho
a}p_{t}\left( \sigma ,t\right) d\sigma \\
&=&-\int_{0}^{\infty }a^{i}e^{-a\rho }p_{a}\left( a,t\right) da-\left( \mu
_{0}+\Psi \left( P\left( t\right) \right) \right) \int_{0}^{\infty
}a^{i}e^{-a\rho }p\left( a,t\right) da \\
&=&\int_{0}^{\infty }\left( ia^{i-1}e^{-a\rho }-a^{i}\rho e^{-a\rho }\right)
p\left( a,t\right) da-\left( \mu _{0}+\Psi \left( P\left( t\right) \right)
\right) \int_{0}^{\infty }a^{i}e^{-a\rho }p\left( a,t\right) da \\
&=&i\int_{0}^{\infty }a^{i-1}e^{-a\rho }p\left( a,t\right)
da-\int_{0}^{\infty }a^{i}\rho e^{-a\rho }p\left( a,t\right) da-\left( \mu
_{0}+\Psi \left( P\left( t\right) \right) \right) \int_{0}^{\infty
}a^{i}e^{-a\rho }p\left( a,t\right) da \\
&=&iP_{i}\left( t\right) -\rho P_{i+1}\left( t\right) -\left( \mu _{0}+\Psi
\left( P\left( t\right) \right) \right) P_{i+1}\left( t\right) \\
&=&iP_{i}\left( t\right) -\left( \rho +\mu _{0}+\Psi \left( P\left( t\right)
\right) \right) P_{i+1}\left( t\right) .
\end{eqnarray*}%
Finally, to achieve our goal of obtaining the existence of solutions to the
model (\ref{SNPM}) we are led to the system of differential equations of
first order 
\begin{equation}
\left\{ 
\begin{array}{l}
P^{\prime }\left( t\right) =-\left( \mu _{0}+\Psi \left( P\left( t\right)
\right) \right) P\left( t\right) +R_{0}\Phi \left( P\left( t\right) \right) 
\overset{n-1}{\underset{i=0}{\Sigma }}\beta _{i}P_{i+1} \\ 
P_{1}^{\prime }\left( t\right) =\left( R_{0}\beta _{0}\Phi \left( P\left(
t\right) \right) -\rho -\mu _{0}-\Psi \left( P\left( t\right) \right)
\right) P_{1}\left( t\right) +R_{0}\Phi \left( P\left( t\right) \right) 
\overset{n-1}{\underset{i=1}{\Sigma }}\beta _{i}P_{i+1} \\ 
P_{i+1}^{\prime }\left( t\right) =iP_{i}\left( t\right) -\left( \rho +\mu
_{0}+\Psi \left( P\left( t\right) \right) \right) P_{i+1}\left( t\right) 
\text{, }i=1,...,n-1%
\end{array}%
\right.  \label{stud}
\end{equation}%
coupled with the initial conditions%
\begin{equation}
P\left( 0\right) =P_{0}\text{, }P_{i}\left( 0\right) =P_{0}^{i},\text{ }%
i=1,2,...,n  \label{init}
\end{equation}%
where%
\begin{equation*}
P\left( 0\right) =\int_{0}^{\infty }p_{0}\left( a\right) da\text{, }%
P_{i}\left( 0\right) =\int_{0}^{\infty }a^{i-1}e^{-\rho a}p_{0}\left(
a\right) da,\text{ }i=1,2,..n.
\end{equation*}

The study of the existence of solutions for the system (\ref{stud}) is
equivalent to the study of existence of solutions to (\ref{SNPM}) because,
if the pair 
\begin{equation*}
(P\left( t\right) ,P_{1}\left( t\right) ,...,P_{n}\left( t\right) )
\end{equation*}
solves (\ref{stud}), then by setting%
\begin{equation*}
B\left( t\right) =p\left( 0,t\right) =R_{0}\Phi \left( P\left( t\right)
\right) \overset{n-1}{\underset{i=0}{\Sigma }}\beta _{i}P_{i+1}\left(
t\right)
\end{equation*}%
we obtain the solution to (\ref{SNPM}) via the usual formula%
\begin{equation*}
p\left( a,t\right) =\left\{ 
\begin{array}{cc}
p_{0}\left( a-t\right) e^{-\int_{0}^{t}\left( \mu _{0}+\Psi \left( P\left(
\sigma \right) \right) \right) d\sigma } & a\geq t, \\ 
B\left( t-a\right) e^{-\int_{t-a}^{t}\left( \mu _{0}+\Psi \left( P\left(
\sigma \right) \right) \right) d\sigma } & a<t.%
\end{array}%
\right.
\end{equation*}%
Let us summarize the results here.

\begin{theorem}
The system (\ref{SNPM}) is equivalent to the ordinary nonlinear system (\ref%
{stud}).
\end{theorem}

Thus, in order to determine the existence of solutions to the model (\ref%
{SNPM}) we can focus on the analysis of (\ref{stud}).

Clearly, the system (\ref{stud}) has at least the trivial solution and so
the existence of stationary solutions to this problem (\ref{stud}) are of
our interest that may not be unique for some values of the parameters and
then may lead to complex bifurcations. In the model of the problem any
stationary solution is called an equilibrium density function.

\subsection{The Equilibrium Solution and their dynamic behavior}

The equilibrium solution $\left( P^{\ast },P_{1}^{\ast },...,P_{n}^{\ast
}\right) $ of (\ref{stud}) is given by%
\begin{equation}
\left\{ 
\begin{array}{l}
0=-\left( \mu _{0}+\Psi \left( P^{\ast }\right) \right) P^{\ast }+R_{0}\Phi
\left( P^{\ast }\right) \overset{n-1}{\underset{i=0}{\Sigma }}\beta
_{i}P_{i+1}^{\ast } \\ 
0=\left( R_{0}\beta _{0}\Phi \left( P^{\ast }\right) -\rho -\mu _{0}-\Psi
\left( P^{\ast }\right) \right) P_{1}^{\ast }+R_{0}\Phi \left( P^{\ast
}\right) \overset{n-1}{\underset{i=1}{\Sigma }}\beta _{i}P_{i+1}^{\ast } \\ 
0=iP_{i}^{\ast }-\left( \rho +\mu _{0}+\Psi \left( P^{\ast }\right) \right)
P_{i+1}^{\ast }\text{, }i=1,...,n-1.%
\end{array}%
\right.  \label{studdd}
\end{equation}%
To solve the nonlinear algebraic system (\ref{studdd}), in our attention is
the third equation ($i=1,...,n-1)$ from where we obtain successively 
\begin{eqnarray*}
P_{i+1}^{\ast } &=&\frac{i}{\rho +\mu _{0}+\Psi \left( P^{\ast }\right) }%
P_{i}^{\ast } \\
&=&\frac{i\left( i-1\right) }{\rho +\mu _{0}+\Psi \left( P^{\ast }\right) }%
P_{i-1}^{\ast } \\
&=& \\
&&... \\
&=&\frac{i!}{\left( \rho +\mu _{0}+\Psi \left( P^{\ast }\right) \right) ^{i}}%
P_{1}^{\ast }.
\end{eqnarray*}%
The next step, is to replace the determined quantities%
\begin{equation}  \label{eu}
P_{i+1}^{\ast }=\frac{i!}{\left( \rho +\mu _{0}+\Psi \left( P^{\ast }\right)
\right) ^{i}}P_{1}^{\ast }\text{, }i=1,...,n-1.
\end{equation}%
in the first equation of (\ref{studdd}). By equivalence, we obtain%
\begin{eqnarray}
0 &=&-\left( \mu _{0}+\Psi \left( P^{\ast }\right) \right) P^{\ast
}+R_{0}\Phi \left( P^{\ast }\right) \beta _{0}P_{1}^{\ast }+R_{0}\Phi \left(
P^{\ast }\right) \overset{n-1}{\underset{i=1}{\Sigma }}\beta
_{i}P_{i+1}^{\ast }  \notag \\
&\Leftrightarrow &  \notag \\
0 &=&-\left( \mu _{0}+\Psi \left( P^{\ast }\right) \right) P^{\ast
}+R_{0}\Phi \left( P^{\ast }\right) \beta _{0}P_{1}^{\ast }+R_{0}\Phi \left(
P^{\ast }\right) \overset{n-1}{\underset{i=1}{\Sigma }}\frac{\beta _{i}i!}{%
\left( \rho +\mu _{0}+\Psi \left( P^{\ast }\right) \right) ^{i}}P_{1}^{\ast }
\notag \\
&\Leftrightarrow &  \notag \\
0 &=&-\left( \mu _{0}+\Psi \left( P^{\ast }\right) \right) P^{\ast }+\left[
\beta _{0}+\overset{n-1}{\underset{i=1}{\Sigma }}\frac{\beta _{i}i!}{\left(
\rho +\mu _{0}+\Psi \left( P^{\ast }\right) \right) ^{i}}\right] R_{0}\Phi
\left( P^{\ast }\right) P_{1}^{\ast }.  \notag
\end{eqnarray}%
Finally, since the non-trivial stationary sizes $P^{\ast }$ is given by (\ref%
{verif}), we obtain from equation above that%
\begin{equation}  \label{eeu}
P_{1}^{\ast }=\frac{\left( \mu _{0}+\Psi \left( P^{\ast }\right) \right) }{%
\left[ \beta _{0}+\overset{n-1}{\underset{i=1}{\Sigma }}\frac{\beta _{i}i!}{%
\left( \rho +\mu _{0}+\Psi \left( P^{\ast }\right) \right) ^{i}}\right]
R_{0}\Phi \left( P^{\ast }\right) }P^{\ast }=\frac{\left( \mu _{0}+\Psi
\left( P^{\ast }\right) \right) }{\frac{\rho +\mu _{0}+\Psi \left( P^{\ast
}\right) }{R_{0}\Phi \left( P^{\ast }\right) }R_{0}\Phi \left( P^{\ast
}\right) }=\frac{\mu _{0}+\Psi \left( P^{\ast }\right) }{\rho +\mu _{0}+\Psi
\left( P^{\ast }\right) }P^{\ast }
\end{equation}%
The existence of a non-trivial stationary solution for the system (\ref{stud}%
), in the form%
\begin{equation}
P_{1}^{\ast }=\frac{\mu _{0}+\Psi \left( P^{\ast }\right) }{\rho +\mu
_{0}+\Psi \left( P^{\ast }\right) }P^{\ast }\text{ and }P_{i+1}^{\ast }=%
\frac{i!}{\left( \rho +\mu _{0}+\Psi \left( P^{\ast }\right) \right) ^{i}}%
P_{1}^{\ast }\text{, }i=1,...,n-1,  \label{nontriv}
\end{equation}%
is proved if the second equation in (\ref{studdd}) is checked by (\ref%
{nontriv}). But, this is a simple exercises by the following equivalence%
\begin{eqnarray*}
0 &=&\left( R_{0}\beta _{0}\Phi \left( P^{\ast }\right) -\rho -\mu _{0}-\Psi
\left( P^{\ast }\right) \right) P_{1}^{\ast }+R_{0}\Phi \left( P^{\ast
}\right) \overset{n-1}{\underset{i=1}{\Sigma }}\beta _{i}P_{i+1}^{\ast } \\
&\Longleftrightarrow & \\
0 &=&\left( R_{0}\beta _{0}\Phi \left( P^{\ast }\right) -\rho -\mu _{0}-\Psi
\left( P^{\ast }\right) \right) P_{1}^{\ast }-R_{0}\Phi \left( P^{\ast
}\right) \beta _{0}P_{1}^{\ast }+R_{0}\Phi \left( P^{\ast }\right) \overset{%
n-1}{\underset{i=0}{\Sigma }}\beta _{i}P_{i+1}^{\ast } \\
&\Leftrightarrow & \\
0 &=&\left( -\rho -\mu _{0}-\Psi \left( P^{\ast }\right) \right) P_{1}^{\ast
}+\left[ \beta _{0}+\overset{n-1}{\underset{i=1}{\Sigma }}\frac{\beta _{i}i!%
}{\left( \rho +\mu _{0}+\Psi \left( P^{\ast }\right) \right) ^{i}}\right]
R_{0}\Phi \left( P^{\ast }\right) P_{1}^{\ast } \\
&\Leftrightarrow & \\
0 &=&\left( -\rho -\mu _{0}-\Psi \left( P^{\ast }\right) \right) P_{1}^{\ast
}+\frac{\rho +\mu _{0}+\Psi \left( P^{\ast }\right) }{R_{0}\Phi \left(
P^{\ast }\right) }R_{0}\Phi \left( P^{\ast }\right) P_{1}^{\ast }
\end{eqnarray*}%
and the last equality is true. Let us summarize our result

\begin{lemma}
The equilibrium solution $\left( P^{\ast },P_{1}^{\ast },...,P_{n}^{\ast
}\right) $ of (\ref{stud}) is given by \eqref{ec}, \eqref{eeu}, and %
\eqref{eu}.
\end{lemma}

Let us point that, the stability of the equilibrium point $(P^{\ast
},P_{1}^{\ast },...,P_{n}^{\ast })$ is determined by the sign of trace of
the Jacobian matrix of the system. Depending on the parameter combinations
chosen, the model can show stability as well as instability of the
non-trivial equilibrium. Also, existence of periodic solutions occurs when
passing from one case to the other.

\begin{remark}
Our model analysis can be extended to more general mortality 
\begin{equation*}
\mu (a,p)=\mu_0 (t)+\Psi (p),
\end{equation*}%
where $\mu_0 (t)$ is a deterministic function, i.e., Gompertz function 
\begin{equation*}
\mu_0 (t)=a+be^{-ct},
\end{equation*}%
for some constants $a,b,$ and $c$. However in such a case the trivial
equilibrium is the only equilibrium.

Let us point out that our approach can be applied to more general fertility function

\begin{equation*}
\beta \left( a,p\right) =R_{0}\beta _{0}\Phi \left( p\right) F(a),
\end{equation*}%

for any continuous function $F(a).$ Indeed, this is the case since $F(a)$
can be approximated by

$$ \overset{n-1}{%
\underset{i=0}{\Sigma }}\beta _{i}a^{i}e^{-\rho a}.$$
\end{remark}

\section*{Acknowledgments}

This research was supported by NSERC grant 5-36700 by Traian A. Pirvu and
Horizon2020-2017-RISE-777911 project by Catalin Sterbeti.

\bigskip

\end{document}